\documentstyle{amsppt}
\TagsOnRight
\magnification=\magstep1
\NoRunningHeads
\footline={\hss \tenrm \folio \hss }

\topmatter
\title
ON THE ISOMORPHISM OF SOME  ONE-RELATOR GROUPS
\endtitle
\author
A.~V.~BORSCHEV, D.~I.~MOLDAVANSKII
\endauthor
\abstract
A class of one-relator groups such that every group in the
class is determined by a triple of integers and is an HNN-extension of
some Baumslag -- Solitar group is considered. A criterion for two groups
in this class to
be isomorphic and criterion for two non-isomorphic groups to be a
homomorphic image of each other are obtained. This result leads to
a negative answer to the question 3.33 from "Kourovka Notebook".
\endabstract
\endtopmatter

\document

In the present paper we consider one-relator groups of form
$$
G(l,m;k)=\langle a, t; \ t^{-1}a^{-k}ta^lt^{-1}a^kt=a^m \rangle,
$$
where $l$, $m$ and $k$ are arbitrary nonzero integers.
Investigation of some properties of these groups was undertaken by
A.~M.~Brunner [1] who observed, in particular, that any such group
is isomorphic to some group $G(l,m;k)$ where integers $l$, $m$ and
$k$ satisfy the conditions $k>0$ and $|l|\geqslant m>0$ (and these
conditions will be assumed in what follows). However, earlier, in 1969,
these groups have appeared in the
G.~Baumslag's paper [2] where it was proved that any finite image of
group $G(2,1;1)$ must be cyclic (and thus the most impressive
example of one-relator group which is not residually finite was given).

Later on the results of Brunner were supplemented in [3]: in the case
$|l|>m$ the generators and defining relations of automorphism group of
group $G(l,m;k)$ were given (in [1] this was done in the cases when either
$m=1$ or $m$ does not divide $l$), it was proved also that the group
$G(l,m;k)$ is residually finite if and only if
$|l|=m$ (in [1] the part \lq\lq only if\rq\rq was singled out) and it is
non-Hopfian if and only if $|l|>m>1$, $m$ divides both integers $l$ and $k$
and integers $m$ and $l/m$ are coprime (in [1] the part \lq\lq if\rq\rq was
proved). Then, in [4] was shown that the group $G(l,m;k)$ is residually
a finite $p$-groups if and only if $|l|=m=p^r$ and $k=p^s$ for some integers
$r\geqslant 0$ and $s\geqslant 0$ and also if $l=-m$ then $p=2$ and
$s\leqslant r$.

In this paper we investigate conditions for isomorphism of groups $G(l,m;k)$.
Our main result is the following

\proclaim{\indent Theorem} {\rm 1)} If the groups $G_1=G(l_1,m_1;k_1)$ and
$G_2=G(l_2,m_2;k_2)$ are homomorphic images of each other and at least
one of inequalities $|l_1|>m_1$ and $|l_2|>m_2$ is satisfied then
$l_1=l_2$ and $m_1=m_2$.

{\rm 2)} If $|l|>m$ then the groups $G_1=G(l,m;k_1)$ and $G_2=G(l,m;k_2)$
are isomorphic if and only if one of the following conditions is
satisfied:
\roster
\item"(2.1)" $k_1=k_2$;
\item"(2.2)" $m>1$, the integers $k_1$ and $k_2$ are divided by the
greatest common divisor of $l$ and $m$ and $k_1/k_2=\pm (l/m)^p$ for some
integer $p\neq 0$;
\item"(2.3)" $m=1$ and the quotient $k_1/k_2$ is an $l$-number.
\endroster

{\rm 3)} The groups $G_1=G(l,m;k_1)$ and $G_2=G(l,m;k_2)$
are homomorphic images of each other and are not isomorphic if and only if
$|l|>m>1$, $m$ divides integers $l$, $k_1$ and $k_2$, integers
$s=l/m$ and $m$ are coprime and the quotient $k_1/k_2$ is an $s$-number
not coinciding with any power (with integer exponent) of $\pm s$.
\endproclaim

(Here, as usually, the integer $m$ is called an $n$-number if any prime
divizor of $m$ divides the integer $n$; the rational number will be called
an $n$-number if numerator and denumerator of irreducible fraction
representing it are $n$-numbers.)
\smallskip

The attempt to give the necessary and sufficient conditions for groups from
class in question to be isomorphic was announced in [5]. Nevertheless, the
theorem given
there appears to be somewhat not exact and the item 2) of our Theorem
can be considered as improving of it in the case when $|l_1|>m_1$ or
$|l_2|>m_2$. The question on the conditions for groups $G_1$ and
$G_2$ to be isomorphic in the remaining case $|l_1|=m_1$ and $|l_2|=m_2$
is still open.

Let us remark, however, that if groups $G_1$ and $G_2$ are homomorphic
images of each other and, say, $|l_1|=m_1$, then group $G_1$, being
residually finite, is Hopfian and therefore group $G_2$ must be
isomorphic to $G_1$. Consequently group $G_2$ is residually finite and
therefore $|l_2|=m_2$. If, in addition, we assume that $l_1=-m_1$ then
factorization by commutator subgroups gives $l_1=l_2$ and $m_1=m_2$.
Thus the statement of item 1) of our Theorem is valid also in the case
when at least for one value of $i=1, 2$ the equation $l_i=-m_i$ is satisfied.
As well we see that the satisfiability of inequality $|l|>m$ in the
statement of item 3) in our Theorem is evident.

On the other hand, the corollary stated in [5] is valid. It gives the
sufficient condition for non-isomorphic groups $G_1$ and $G_2$ to be
homomorphic images of each other (necessary and sufficient condition is
stated in item 3) of our Theorem). This leads, in turn, to the negative
answer to the
question 3.33 from [6]: Are two groups necessarily isomorphic if each
of them may be defined by a single relation and is a homomorphic image
of the other one? In fact, even the condition in [5]
implies that groups $G(18,2;2)$ and $G(18,2;6)$ give counterexample, while
the item 3) of our Theorem implies that this counterexample is in some
sense minimal.
\smallskip

The proof of Theorem occupies the section 2 of the paper. Section 1
contains some required for the proof results that either are contained
in papers [1] and [3] or are immediate consequences of them.
We emphasize once more that everywhere below it is assumed without additional
reminder that the integers $l$, $m$ and $k$ (with the same subscripts or
without subscripts) do satisfy the inequalities $|l|>m>0$ and $k>0$
(the condition $|l|>m>0$ is essential, in particular, in references to
results in [1]).
\bigskip

\centerline{{\bf 1. Preliminaries}}
\medskip

Following [1], add to the presentation
$\langle a, t; \ t^{-1}a^{-k}ta^lt^{-1}a^kt=a^m \rangle$
of group $G(l,m;k)$ new generator $b$ and defining relation
$b=t^{-1}a^kt$. Resulting presentations
$$
G(l,m;k)=\langle a, b, t; \ b^{-1}a^lb=a^m, t^{-1}a^kt=b \rangle
$$
obviously shows that the group $G(l,m;k)$ is an $HNN$-extension
with stable letter $t$ of base group
$H(l,m)=\langle a, b; \ b^{-1}a^lb=a^m \rangle$, where associated
subgroups $H^{-1}$ and $H^1$ are infinite cyclic with generators
$a^k$ and $b$ respectively (all relevant definitions and assertions
relative to construction of $HNN$-extension can be found in [7]).
Any group $H(l,m)$ (known as Baumslag -- Solitar group)
is, in turn, $HNN$-extension with stable letter $b$ of infinite
cyclic group generated by element $a$. Therefore it will be necessary
to speak about $t$-length and  $t$-reduced form and about
$b$-length and  $b$-reduced form of elements of group $G(l,m;k)$.

Symbol $L$ will denote the normal closure in group $H(l,m)$ of element $a$.
It is obvious that the group $H(l,m)$ is a split extension of subgroup
$L$ by infinite cyclic group generated by $b$ and that element
$h\in H(l,m)$ belongs to $L$ if and only if the exponent sum on $b$ in
(any) word in $a$ and $b$ representing $h$ is equal to $0$.
Let also $C(h)$ denote the centralizer in group $H(l,m)$ of element $h$.

When two groups of form $G(l,m;k)$ are considered simultaneously we will
write $G_1=G(l_1,m_1;k_1)$ and $G_2=G(l_2,m_2;k_2)$ and will equip by
corresponding subscripts generators $a, b, t$ and subgroups $H=H(l,m)$,
$H^{-1}$, $H^1$ and $L$.

Let us state now certain relevant properties of groups $H(l,m)$ and $G(l,m;k)$.
The immediate induction on $b$-length of element $h$ gives the following
assertion which is a some refinement of Lemma 2.4 from [1]:

\proclaim{\indent Proposition 1} Let element $h\in H(l,m)$ be written
in the form $h=b^pv$, where $v\in L$, $p\in \Bbb Z$. Let $d=(l,m)$ be
the greatest common divizor of integers $l$ and $m$ and $l=l_1d$, $m=m_1d$.
The equality $h^{-1}a^rh=a^s$ is satisfied if and only if
$v\in C(a^s)$ and either $r=s$ if $p=0$, or $r=l_1^pdx$ and $s=m_1^pdx$
for some integer $x$ if $p>0$, or $r=m_1^{-p}dx$ and $s=l_1^{-p}dx$
for some integer $x$ if $p<0$.
\endproclaim

\proclaim{\indent Proposition 2} {\rm (Lemma 2.3 from [1])} Let
$h\in H(l,m)$. If for some integer $n\neq 0$ element $h^n$ is conjugate to
element from subgroup generated by $a$ or from subgroup generated by $b$
then element $h$ itself is conjugate to an element from the same subgroup.
\endproclaim

\proclaim{\indent Proposition 3} Let for some element $g\in G(l,m;k)$
and for some integer $n\neq 0$ the inclusion $g^{-1}a^ng\in H(l,m)$
be fulfilled. Then element $g^{-1}a^ng$ belongs to subgroup $L$ if and
only if element $g$ belongs to subgroup $H(l,m)$. In addition, if
$g\notin H(l,m)$ then elements $g$ and $g^{-1}a^ng$ are of form
$g=b^putv$ and $g^{-1}a^ng=v^{-1}b^sv$, where  $s\neq 0$ and $p$ are
some integers and $u$ and $v$ are some elements from  $L$ such that
$u\in C(a^{ks})$ and $b^{-p}a^nb^p=a^{ks}$.
\endproclaim

Indeed, Lemma 2.5 of paper [1] implies that if element $g$ such that
$g^{-1}a^ng\in H(l,m)$ for some integer $n\neq 0$ does not belong to
subgroup $H(l,m)$, then it is of form $g=wtv$ where $w\in H(l,m)$ and
$v\in L$. Since the expression $v^{-1}t^{-1}w^{-1}a^nwtv$ of element
$g^{-1}a^ng$ can not be $t$-reduced, we must have $w^{-1}a^nw=a^{ks}$
for some integer $s\neq 0$. Hence $g^{-1}a^ng=v^{-1}b^sv$ and since
$s\neq 0$ this element does not belong to subgroup $L$. If now
element $w$ is represented in the form $w=b^pu$, where $u\in L$, then
by Proposition 1 we have $u\in C(a^{ks})$ and therefore
$b^{-p}a^nb^p=a^{ks}$.
\smallskip

An endomorphism $\varphi$ of group $G(l,m;k)$ we shall call
{\it special} if $a\varphi =a^r$ and $b\varphi=bg$ for some integer
$r\neq 0$ and element $g\in L$. The number $r$ will be named
{\it exponent} of endomorphism $\varphi$. Lemma 3.1 in [1] asserts,
in particular, that for any endomorphism $\psi$ of $G(l,m;k)$ whose image
is not cyclic there exists an inner automorphism $\tau$ of this group
such that the product $\psi\tau$ is a special endomorphism. Furthermore,
certain part of the proof in [1] of this Lemma gives the following

\proclaim{\indent Proposition 4} Let $\varphi$ be a special endomorphism
of the group $G(l,m;k)$ with exponent $r$. Then for some integer
$p\geqslant 0$ the equality $rm^p=l^p$ holds and element $t\varphi$ is of
form $b^putv$ where  $u, v\in L$. Furthermore, $b\varphi=v^{-1}bv$
(and therefore element $g$ from the definition above of special endomorphism
is equal to commutator $[b,v]$).
\endproclaim

\proclaim{\indent Proposition 5}
{\rm 1)} Every special endomorphism of group $G(l,m;k)$ with exponent
$r=1$ is injective.

{\rm 2)} If the group $G(l,m;k)$  has a special endomorphism with
exponent $r\neq \pm 1$ then $m$ divides both integers $l$ and $k$.

{\rm 3)} If $m>1$ and a special endomorphism with exponent $r\neq \pm 1$
is surjective then it is not injective.
\endproclaim

Assertion 1) of this Proposition is contained in Proposition 3.2 from
paper [3] and assertions 2) and 3) coincide respectively with
Propositions 3.5 and 3.6 from the same paper.
\bigskip

\centerline{{\bf 2. The proof of Theorem}}
\medskip

Let the groups $G_1=G(l_1,m_1,k_1)$ and $G_2=G(l_2,m_2,k_2)$ be homomorphic
images of each other and let $\varphi$ be a surjective homomorphism
of $G_1$ onto $G_2$. Then in group $G_2$ must be satisfied the equality
$$
(b_1\varphi)^{-1}(a_1\varphi)^{l_1}(b_1\varphi)=(a_1\varphi)^{m_1}. \tag1
$$
Since $|l_1|>m_1$, Collins' Lemma ([7], p.~185) implies that element
$a_1\varphi$ must be conjugate to some element of subgroup $H_2$.
Therefore  multiplying the mapping $\varphi$ by suitable inner
automorphism of $G_2$ permits to assume without loss of generality that
$a_1\varphi \in H_2$. We want show now that, in fact, it may be even
assumed (again after multiplying $\varphi$ by suitable inner
automorphism of $G_2$) that for some integer $r_1\neq 0$ the equality
$$
a_1\varphi=a_2^{r_1} \tag2
$$
holds in group $G_2$ and that, moreover, $b_1\varphi\in H_2$.

Indeed,
if for mapping $\varphi$ the inclusion $b_1\varphi\in H_2$ is already
satisfied then the equality (1) and Collins' Lemma (applied to the group
$H_2$) imply that element $a_1\varphi$ is conjugate in $H_2$ to some
power of element $a_2$.

But if $t_2$-length of element $b_1\varphi$ is positive then the equality
(1) and Britton's Lemma imply that element $(a_1\varphi)^{l_1}$ is
conjugate in $H_2$ to some element of subgroup $H_2^{-1}$ or of subgroup
$H_2^1$. By Proposition 2 then element $a_1\varphi$ is conjugate to some
power of element $a_2$ or to some element of subgroup $H_2^1$.
Since subgroups $H_2^{-1}$ and $H_2^1$ are conjugate in group $G_2$,
the existence of integer $r_1$ satisfying the equality (2) is proved in
this case too. In addition, since the mapping
$\varphi$ is surjective and group $G_2$ is non-cyclic, we have $r_1\neq 0$.
Finally, since by (1) we have now the inclusion
$(b_1\varphi)^{-1}a_2^{r_1l_1}(b_1\varphi)\in L_2$, Proposition 3
implies that $b_1\varphi\in H_2$.

Now since elements $a_1\varphi$ and $b_1\varphi$ are in subgroup $H_2$
and homomorphism $\varphi$ is surjective, element $t_1\varphi$ is not
contained in $H_2$. Furthermore, the relation
$$
(t_1\varphi)^{-1}(a_1\varphi)^{k_1}(t_1\varphi)=b_1\varphi
$$
satisfied in group $G_2$ now takes the form
$(t_1\varphi)^{-1}a_2^{r_1k_1}(t_1\varphi)=b_1\varphi$ and therefore
Proposition 3 implies that elements
$t_1\varphi$ and
$(t_1\varphi)^{-1}a_2^{r_1k_1}(t_1\varphi)$ are of form
$b_2^{p_1}u_1t_2v_1$ and $v_1^{-1}b_2^{s_1}v_1$ respectively
for some integers $s_1\neq 0$ and  $p_1$ and for some elements
$u_1$ and $v_1$ from subgroup $L_2$ such that $u_1\in C(a_2^{k_2s_1})$  and
$$
b_2^{-p_1}a_2^{r_1k_1}b_2^{p_1}=a_2^{k_2s_1}. \tag3
$$
Consequently $b_1\varphi=v_1^{-1}b_2^{s_1}v_1$ and relation (1) takes
the form
$$
v_1^{-1}b_2^{-s_1}v_1a_2^{r_1l_1}v_1^{-1}b_2^{s_1}v_1=a_2^{r_1m_1}. \tag4
$$
Now Proposition 1
(taking into account inequalities $|l_i|>m_i$ ($i=1,2$)) gives that
$s_1>0$ and for suitable integer $x$ equalities
$$
r_1l_1=l_{21}^{s_1}d_2x, \qquad r_1m_1=m_{21}^{s_1}d_2x, \tag5
$$
are satisfied where $d_2=(l_2,m_2)$, $l_2=l_{21}d_2$, $m_2=m_{21}d_2$.
Hence we have
$$
l_1/m_1=(l_2/m_2)^{s_1}. \tag6
$$

The same arguments can be applied to any surjective homomorphism
$\psi$ of group $G_2$ onto $G_1$. In particular, we can assume that
for some integers $r_2\neq 0$ and $s_2>0$ and element $v_2\in L_1$ the
equalities $a_2\psi=a_1^{r_2}$ and $b_2\psi=v_2^{-1}b_1^{s_2}v_2$ are
satisfied and that
$$
l_2/m_2=(l_1/m_1)^{s_2}. \tag7
$$

Equalities (6) and (7) implies that $l_1/m_1=(l_1/m_1)^{s_1s_2}$ and
since $|l_1|>m_1$ we obtain $s_1=s_2=1$ and $l_1/m_1=l_2/m_2$.

As a result of these arguments we have the

\proclaim{\indent Lemma 1} Let groups $G_1=G(l_1,m_1,k_1)$ and
$G_2=G(l_2,m_2,k_2)$ be homomorphic images of each other. Then
$l_1/m_1=l_2/m_2$ and up to multipliers which are inner automorphisms
of corresponding groups arbitrary surjective homomorphisms
$\varphi$ of group $G_1$ onto $G_2$ and $\psi$ of group $G_2$ onto $G_1$
are such that $a_1\varphi=a_2^{r_1}$,
$b_1\varphi=v_1^{-1}b_2v_1$ and $a_2\psi=a_1^{r_2}$,
$b_2\psi=v_2^{-1}b_1v_2$ for some nonzero integers
$r_1$ and $r_2$ and for some elements $v_1\in L_2$ and $v_2\in L_1$.
In addition, for some integers $p_1$ and $p_2$ relations
$$
\gather
b_2^{-p_1}a_2^{r_1k_1}b_2^{p_1}=a_2^{k_2}  \quad \text{and} \quad
b_1^{-p_2}a_1^{r_2k_2}b_1^{p_2}=a_1^{k_1}  \tag8 \\
v_1^{-1}b^{-1}_2v_1a_2^{r_1l_1}v_1^{-1}b_2v_1=a_2^{r_1m_1}
\quad \text{and} \quad
v_2^{-1}b^{-1}_1v_2a_1^{r_2l_2}v_2^{-1}b_1v_2=a_1^{r_2m_2} \tag9
\endgather
$$
are satisfied and for some integers $x$, $y$  equalities
$$
r_1l_1=l_2x, \quad r_1m_1=m_2x \quad \text{and} \quad r_2l_2=l_1y,
\quad r_2m_2=m_1y \tag10
$$
hold. Each of endomorphisms $\varphi\psi$ of group $G_1$ and $\psi\varphi$ of
group $G_2$ is special with exponent $r_1r_2$ and therefore for some
integer
$p\geqslant 0$ equalities
$$
r_1r_2m_1^p=l_1^p \quad \text{and} \quad r_1r_2m_2^p=l_2^p. \tag11
$$
are satisfied.
\endproclaim

Since equalities (8), (9) and (10) follows from equalities
(3), (4) and (5) (and from their analogs for mapping $\psi$), we must
prove only the last assertion of Lemma 1. Let us show, for example, that
the mapping $\varphi\psi$ is a special endomorphism of group
$G_1$ with exponent $r_1r_2$. Indeed, the equality
$a_1(\varphi\psi)=a_1^{r_1r_2}$ is evident and since
$b_1\varphi=v_1^{-1}b_2v_1$ and $b_2\psi=v_2^{-1}b_1v_2$ where
$v_1\in L_2$ and $v_2\in L_1$, we have
$$
b_1(\varphi\psi)=(v_1^{-1}b_2v_1)\psi=(v_2(v_1\psi))^{-1}b_1(v_2(v_1\psi))
=b_1g,
$$
where $g=[b_1,v_2(v_1\psi)]\in L_1$. The existence of integer
$p\geqslant 0$ satisfying the equalities (11) follows now from
Proposition 4 and from equality $l_1/m_1=l_2/m_2$.
\smallskip

In the following three Lemmas will be assumed
that $\varphi$ and $\psi$ are surjective homomorphism
of group $G_1$ onto $G_2$ and of group $G_2$ onto $G_1$ respectively
having form indicated in Lemma 1 and also will be used corresponding
notations and relations from this Lemma.

\proclaim{\indent Lemma 2} Let $|r_1r_2|=1$. Then $l_1=l_2$, $m_1=m_2$
and if $k_1\neq k_2$, then integers $k_1$ and $k_2$ are divided by
the greatest common divizor of
$l_1$ and $m_1$ and $k_1/k_2=\pm (l_1/m_1)^{q}$ for some integer
$q\neq 0$.
\endproclaim

Indeed, since $r_1=\pm 1$, from equalities (10) it follows, in particular,
that $\pm m_1=m_2x$ and therefore $m_2$ divides $m_1$. Similarly,
$m_1$ divides $m_2$, and so $m_1=m_2$. Since $l_1/m_1=l_2/m_2$,
we have $l_1=l_2$. In addition, the first from relations (8) gives
the equality
$b_2^{-p_1}a_2^{\varepsilon k_1}b_2^{p_1}=a_2^{k_2}$ (where
$\varepsilon=\pm 1$), which by Proposition 1 (and with regard to
positiveness of numbers $k_1$ and $k_2$) is valid if and only if
either $k_1=k_2$, or $p_1\neq 0$, numbers $k_1$ and $k_2$ are divided by
the greatest common divizor of $l_1$ and $m_1$ and
$k_1/k_2=\varepsilon (l_1/m_1)^{p_1}$.
\smallskip

\proclaim{\indent Lemma 3} Let $|r_1r_2|>1$, $m_1>1$ and $m_2>1$.
Then $l_1=l_2$, $m_1=m_2$, $m_1$ divides numbers
$l_1$, $k_1$ and $k_2$, integer $s=l_1/m_1$ is coprime to
$m_1$ and quotient $k_1/k_2$ is $s$-number.
\endproclaim

\demo{\indent Proof} Since endomorphism $\varphi\psi$ of group
$G_1$ is special with exponent $r_1r_2$, the item 2) of Proposition 5
imlies that $m_1$ must divide numbers
$l_1$ and $k_1$. If $l_1=m_1s$, then by (11) we have
$r_1r_2=s^p$ where $p>0$. We claim that integers $m_1$ and $s$
are coprime.

Since endomorphism $\varphi\psi$ is surjective,
for some element $w\in G_1$ we must have $a_1=w(\varphi\psi)$.
If element $w$ is written as a word
$w(a_1, b_1, t_1)$ in generators of group $G_1$ we have the equality
$a_1=v(a_1, b_1, t_1)$ where the word
$$
v(a_1, b_1, t_1)=w(a_1^{s^p}, b_1(\varphi\psi), t_1(\varphi\psi))
$$
is obtained from word $w(a_1, b_1, t_1)$ by replacement of every
generator by some word that is equal to
$\varphi\psi$-image of this generator.

Let us calculate the exponent sum on $a_1$ in word $v(a_1, b_1, t_1)$.
To do this we note, first of all, that by Proposition 4
we may assume that in the word representing element $b_1(\varphi\psi)$
exponent sum on $a_1$ is equal to zero. Next, the form of element
$t_1(\varphi\psi)$ stated in Proposition 4 shows that
exponent sum on $t_1$ in word $v(a_1, b_1, t_1)$ coincide with exponent
sum on $t_1$ in word $w(a_1, b_1, t_1)$. On the other hand,
the exponent sum on $t_1$ in any word which is equal to identity in group
$G_1$ must be zero. Therefore, the equality
$a_1=v(a_1, b_1, t_1)$ implies that
exponent sum on $t_1$ in word  $v(a_1, b_1, t_1)$ and therefore
in word $w(a_1, b_1, t_1)$ is equal to zero. Thus, the exponent sum on
$a_1$ in word $v(a_1, b_1, t_1)$ is equal to $ns^p$, where $n$ is the
exponent sum on $a_1$ in word $w(a_1, b_1, t_1)$.

Since $m_1$ divides integers $l_1$ and $k_1$, the exponent sum
on $a_1$ in defining relators of group $G_1$ and therefore in any word that
equal to 1 is divided by $m_1$. Consequently, the equality
$a_1=v(a_1, b_1, t_1)$ implies the congruence $ns^p\equiv 1\pmod {m_1}$
which, in turn, implies that integers $m_1$ and $s$ are coprime.

The same arguments in connection with endomorphism
$\psi\varphi$ of group $G_2$ (with regard to equality $l_1/m_1=l_2/m_2$)
give the equality $l_2=m_2s$ and coprimeness of integers $m_2$ and $s$.

The second from equalities (10) says that $m_2$ divides integer
$r_1m_1$. But since $r_1r_2=s^p$ the integer $r_1$ is an $s$-number
and therefore is coprime to $m_2$. Consequently,  $m_2$ divides
$m_1$ and similarly  $m_1$ divides $m_2$, hence we obtain
$m_1=m_2$ and $l_1=l_2$.

Finally, since $l_1=m_1s$, the first from relations (8) and Proposition 1
imply $r_1k_1=s^{p_1}k_2$. Therefore, $k_1/k_2$ is an $s$-number and
Lemma 3 is proved.

\proclaim{\indent Lemma 4} If one of the numbers $m_1$ and $m_2$ is
equal to 1 then other is equal to 1 too. If $m_1=m_2=1$ then $l_1=l_2$
and $k_1/k_2$ an $l_1$-number.
\endproclaim

Suppose, on the contrary,  that $m_1>1$ and $m_2=1$. Then by Lemma 2
$|r_1r_2|>1$ and therefore by Proposition 5 endomorphism $\varphi\psi$
of group $G_1$ is not injective.

On the other hand from (11) it follows that $r_1r_2=l_2^p$. Therefore
the product of endomorphism $\psi\varphi$ and of inner automorphism
of group $G_2$ generated by $b_2^p$ fixed element $a_2$,
i.~e. this product is special endomorphism with exponent 1.
Proposition 5 now implies that endomorphism $\psi\varphi$ is injective
and hence with regard to surjectiveness $\psi$ it follows injectiveness
of mappings $\varphi$ and $\psi$. So, the mapping $\varphi\psi$ is
injective and this contradiction proves the first assertion of Lemma.

If now $m_1=m_2=1$ then the equality $l_1=l_2$ is evident. Since in that
case the first of relations (8) is equivalent to equality
$r_1k_1=l_1^{p_1}k_2$ and by (11) $r_1$ is $l_1$-number, the quotient
$k_1/k_2$ is $l_1$-number and Lemma 4 is proved.
\enddemo

Now we can proceed directly to prove all assertions of Theorem.

Let the groups $G_1=G(l_1,m_1,k_1)$ and $G_2=G(l_2,m_2,k_2)$ be homomorphic
images of each other. Then by Lemma 1 $l_1/m_1=l_2/m_2$ and hence if
$m_1=m_2=1$ then  $l_1=l_2$. Otherwise, by Lemma 4 $m_1>1$ and $m_2>1$,
and in this case assertion of item 1) follows from Lemmas 2 and 3.

To prove assertion of item  2) let us at first suppose that groups
$G_1=G(l,m,k_1)$ and $G_2=G(l,m,k_2)$ are isomorphic. Let $\varphi$ be
isomorphism of group $G_1$ onto $G_2$ and  $\psi$ be isomorphism of group
$G_2$ onto $G_1$. Without loss of generality it may be supposed
that these mappings are of form indicated by Lemma 1.
Let $k_1\neq k_2$. If $m>1$ then since the mapping
$\varphi\psi$ is bejective it follows from Proposition 5 that $|r_1r_2|=1$.
So, in this case the validity of assertion
(2.2) follows from Lemma 2. The assertion
(2.3) is contained in Lemma 4.

Conversely, let us suppose at first that the condition
(2.2) is satisfied i.~e. integers $k_1$ and $k_2$ are divided by the
greatest common divizor $d$ of $l$ and $m$ and for some
$\varepsilon=\pm 1$ and for some integer $p\neq 0$ the equation
$k_1/k_2=\varepsilon (l/m)^p$ holds. It is easily to see that
for some integer $x$ either $\varepsilon k_1=l_1^pdx$ and $k_2=m_1^pdx$
if $p>0$ or $\varepsilon k_1=m_1^{-p}dx$ and $k_2=l_1^{-p}dx$ if $p<0$
(where $l_1=l/d$ and $m_1=m/d$). It can be immediately verified
(see, however, Proposition 1), that in any case in groups
$G_2$ and $G_1$ the equalities
$b_2^{-p}a_2^{\varepsilon k_1}b_2^p=a_2^{k_2}$ and
$b_1^{p}a_1^{\varepsilon k_2}b_1^{-p}=a_1^{k_1}$ are satisfied.
Therefore the mappings
$$
a_1\mapsto a_2^{\varepsilon},\ b_1\mapsto b_2,\ t_1\mapsto b_2^pt_2
\quad \text{and} \quad
a_2\mapsto a_1^{\varepsilon},\ b_2\mapsto b_1,\ t_2\mapsto b_1^{-p}t_1
$$
define mutually inverse isomorphisms
of group $G_1$ onto $G_2$ and of group $G_2$ onto $G_1$ respectively.

Let now the condition (2.3) be satisfied. Then for suitable integer
$l$-numbers $x$ and $y$ the equality $xk_1=yk_2$ holds.
Let us choose an integer $p>0$ such that $l^p=yz$ for some
integer $z$. Then for $r=xz$ we have $rk_1=l^pk_2$ and immediate verification
shows that in group
$G_2$ the equality $b_2^{-p}a_2^{rk_1}b_2^p=a_2^{k_2}$ is satisfied.
Therefore the mapping
$a_1\mapsto a_2^r,\ b_1\mapsto b_2,\ t_1\mapsto b_2^pt_2$
defines homomorphism of group $G_1$ in group $G_2$. It is easily to see
that since $r$ is $l$-number this homomorphism is surjective. Similarly,
there exists a surjective homomorphism of group
$G_2$ onto $G_1$. Since by [3] these groups are Hopfian they are
isomorphic. The proof of item 2) of Theorem is complete.

To prove the assertion of item 3) let us suppose at first that
groups $G_1=G(l,m,k_1)$ and $G_2=G(l,m,k_2)$ are non-isomorphic and
that $\varphi:G_1\to G_2$ and $\psi:G_2\to G_1$ are surjective homomorphisms
of the form indicated in Lemma 1.
From Lemmas 2 and 4 and from assertions of item 2)proved above  it follows
that $|r_1r_2|>1$ and $m>1$. Therefore the necessity of
conditions in item 3) it follows from Lemma 3 and from assertion of
item (2.2).

Conversely, let us show that if $l=ms$, integers $m$ and $s$ are
coprime, $m$ divides both integers $k_1$ and $k_2$ and
$xk_1=yk_2$ for some integer $s$-numbers $x$ and $y$ then
groups $G_1$ and $G_2$ are homomorphic image of each other.
Let the integer $p>0$ be such that $s^p=yz$ for some integer $z$.
Since then $k_1r=s^pk_2$, where $r=xz$, and $k_1$ and $k_2$ are divided
by $m$, in group $G_2$ the relation
$b_2^{-p}a_2^{k_1r}b_2^p=a_2^{k_2}$ holds. Consequently the mapping
$$
a_1\mapsto a_2^{r}, \ b_1\mapsto b_2, \ t_1\mapsto b_2^pt_2
$$
defines homomorphism of group $G_1$ to group $G_2$. To see that this
homomorphism is surjective it is enough to show that for any
$s$-number $r$ elements $a_2^r$ and $b_2$ generate group $H_2$.
Indeed, let
$n$ be the smallest positive $s$-number such that $a^n$ is contained in
subgroup of $H_2$ generated by these elements. Suppose that
$n>1$ and let
$q$ be the prime divizor of $n$. Then writing $n=n_1q$ and $s=s_1q$ we
see that this subgroup contains element
$b_2^{-1}a_2^{nms_1}b_2=b_2^{-1}a_2^{n_1ms}b_2=a_2^{n_1m}$
and therefore it contains element $a_2^{n_1}$. This contradiction means
that $n=1$ and our claim is proved. Similarly can be proved the existence
of surjective homomorphism of group $G_2$ onto group $G_1$.

It remains to note that if, in addition, $m>1$ and $k_1/k_2$ is not a
power of $\pm s$, then by assertion (2.2) groups $G_1$ and $G_2$ are
not isomorphic.

The proof of Theorem is complete.
\bigskip

\centerline{{\bf References}}
\medskip

\noindent
[1] Brunner A.~M., {\it On a class of one--relator groups.} Can. J. Math.
1980. V. 32. N 2. P. 414--420.

\noindent
[2] Baumslag G., {\it A noncyclic one-relator group all of whose finite
quotients are cyclic.}  J. Austral. Math. Soc. 1969. V. 10. N 3--4.
P. 497--498.

\noindent
[3]  Kavutskii M.~A., Moldavanskii D.~I., {\it On certain class of
one-relator groups.} in \lq\lq Algebraic and Discrete Systems.\rq\rq
Ivanovo State University. Ivanovo, 1988. p. 35--48. (Russian)

\noindent
[4]  Moldavanskii D.~I., {\it Residuality a finite $p$-groups of
$HNN$-extensions.} Vestnik of Ivanovo State University.
2000, Issue 3. p. 129--140. (Russian)

\noindent
[5]  Borschev A.~V. {\it On the isomorphism problem for certain class
of one-relator groups.} International Algebraic conference dedicated to
the memory of D.~K.~Faddeev (St.~Petersburg, Russia, 1997, June 24--30).
Abstracts. St.~Petersburg, 1997. p. 170--171. (Russian)

\noindent
[6] {\it Kourovka Notebook. Unsolved questions in group theory.}
15-th ed. Novosibirsk, 2002.

\noindent
[7] Lindon R., Schupp P. {\it Combinatorial Group Theory.} Springer-Verlag,
Berlin, Heidelberg, New York, 1977.

\end